\documentclass{article}%
\usepackage{amsfonts}
\usepackage{amsmath}
\usepackage{amssymb}
\usepackage{graphicx}%
\newtheorem{theorem}{Theorem}

\newtheorem{definition}[theorem]{Definition}

\newtheorem{proposition}[theorem]{Proposition}
\newtheorem{remark}[theorem]{Remark}

\begin{document}

\title{Explicit solutions of generalized Cauchy-Riemann systems using the transplant operator}
\author{Vladislav V. Kravchenko$^{1}$, S\'{e}bastien Tremblay$^{2}$\\$^{1}$Departamento de Matem\'{a}ticas, CINVESTAV del IPN, Unidad\\Quer\'{e}taro, Libramiento Norponiente No.~2000 C.P. 76230 Fracc.\\Real de Juriquilla, Quer\'{e}taro, Mexico\\$^{2}$D\'{e}partement de math\'{e}matiques et d'informatique, Universit\'{e} du\\Qu\'{e}bec, Trois-Rivi\`{e}res, Qu\'{e}bec, G9A 5H7, Canada}
\maketitle

\begin{abstract}
In \cite{KrTransplant} it was shown that the tool introduced there and called
the transplant operator transforms solutions of one Vekua equation into
solutions of another Vekua equation, related to the first via a
Schr\"{o}dinger equation. In this paper we prove a fundamental property of
this operator: it preserves the order of zeros and poles of generalized
analytic functions and transforms formal powers of the first Vekua equation
into formal powers of the same order for the second Vekua equation. This
property allows us to obtain positive formal powers and a generating sequence
of a \textquotedblleft complicated\textquotedblright\ Vekua equation from
positive formal powers and a generating sequence of a \textquotedblleft
simpler\textquotedblright\ Vekua equation. Similar results are obtained
regarding the construction of Cauchy kernels. Elliptic and hyperbolic
pseudoanalytic function theories are considered and examples are given to
illustrate the procedure.

\end{abstract}

\section{Introduction}

In the present work a special class of Vekua equations describing generalized
analytic or pseudoanalytic functions is considered. It arises naturally in
relation with some linear equations of mathematical physics such as the
stationary Schr\"{o}dinger equation, the conductivity equation and others.
Vekua equations of this type we call main Vekua equations. They are closely
related to another generalization of the Cauchy-Riemann system, the system
describing so-called $p$-analytic functions (see the definitions in the next
section). The general pseudoanalytic function theory mainly created by L. Bers
and his coauthors and presented in \cite{Berskniga} among other developments
contains deep results on generalizations of the concept of complex
differentiability and integrability, Taylor and Laurent series related to
generalized analytic functions as well as the generalizations of the Cauchy
integral formula and its corollaries. In the core of Bers' theory there is a
concept of a generating sequence related to a Vekua equation. In general a
derivative of a generalized analytic function in the sense introduced by Bers
is not any more a solution of the same Vekua equation but of another Vekua
equation called a successor of the original one. Bers derivatives of solutions
of this second Vekua equation will solve another Vekua equation, and in
principle this sequence of Vekua equations related to the original one is
infinite. If somehow one manages to obtain a pair of solutions in a certain
sense independent for each of these Vekua equations then such sequence of
pairs is called the generating sequence and it immediately allows one to
construct a complete system of positive formal powers related to the original
Vekua equation. The formal powers are basic constituents of the Taylor-type
series expansions of the pseudoanalytic functions and generalize the usual
powers $(z-z_{0})^{n}$ in the sense that being a solution of the Vekua
equation a formal power of order $n$ asymptotically behaves like
$(z-z_{0})^{n}$ when $z\rightarrow z_{0}$. Moreover, theorems generalizing
such facts like the Runge theorem on the completeness of the system of
powers\ in a uniform convergence topology and even stronger results
guaranteeing the completeness in the $C$-norm were obtained in the framework
of pseudoanalytic function theory.

One of the most significant obstacles for the further development and a
broader application of pseudoanalytic function theory is the explicit
construction of generating sequences, formal powers and Cauchy kernels
corresponding to Vekua equations arising in applications. Bers himself and
together with Gelbart succeeded in constructing a generating sequence in a
very special case (see \cite{Berskniga} and \cite{APFT}). In
\cite{KrRecentDevelopments} an algorithm for explicit construction of a
generating sequence was proposed for a much more general situation. In
application to second-order elliptic equations with the aid of the tools from
pseudoanalytic function theory this result allows one to obtain a complete
system of solutions of the equation, e.g., of the Schr\"{o}dinger equation
when the equation possesses a particular solution $f$ in a separable form
$f=U(u)V(v)$ where $u,v$ are orthogonal coordinates and $U$, $V$ are arbitrary
sufficiently smooth functions.

In the present paper we substantially extend the class of Vekua equations and
of systems desribing $p$-analytic functions for which a generating sequence
and a system of formal powers can be constructed explicitly. For this we use a
concept introduced in \cite{KrTransplant} and called there the transplant
operator. In fact, it is an operator transforming solutions of one Vekua
equation into solutions of another one related to the first via a
Schr\"{o}dinger equation. Here we prove a crucial property of the transplant
operator: it transforms formal powers into formal powers of the same order
(see details in Section \ref{SectTransplant}). This means that if we are able
to solve a Vekua equation, that is we know its generating sequence then using
the transplant operator we can construct positive formal powers and a
generating sequence for a related Vekua equation which can be much more
complicated. As an example in Section \ref{SectPositiveFormalPowers} we
consider a Vekua equation whose generating sequence is periodic with a period
1, that is it consists of one generating pair only. In this case it is
relatively easy to obtain the corresponding positive formal powers. Then using
the approach described in the present work, it is possible to obtain systems
of positive formal powers and generating sequences for a wide class of Vekua
equations related to the first one. The structure of generating sequences and
of formal powers for the related Vekua equations are more complicated. We also
obtain a similar result regarding the construction of Cauchy kernels as we
show in Section \ref{SectCauchykernels}. If a Cauchy kernel for a Vekua
equation is known, it can be used for constructing Cauchy kernels for a wide
class of related Vekua equations using the transplant operator.

All the described results have a direct application to linear second-order
equations. For example, in the case of the stationary two-dimensional
Schr\"{o}dinger equation $\left(  -\Delta+q(x,y)\right)  u=0$ with $q$ and $u$
being real valued, the existence of one solution $u$ such that a generating
sequence for an associated main Vekua equation can be constructed explicitly
leads not only to the construction of a complete system of solutions to this
Schr\"{o}dinger equation but also to the construction of complete systems of
solutions to any Schr\"{o}dinger equation with the potential $q_{f}%
=-q+2\left(  \nabla f/f\right)  ^{2}$ where $f$ is any solution of the
original Schr\"{o}dinger equation with the potential $q$. Note that the form
of the potential $q_{f}$ is a precise generalization of the potential obtained
after a Darboux transformation in a one-dimensional case (see, e.g., \cite{MS}).

\section{Some known facts about generalized Cauchy-Riemann systems}

Let $\Omega$ be a domain in $\mathbf{R}^{2}$. Throughout the whole paper we
suppose that $\Omega$ is a simply connected domain.

In the present work we consider two related generalized Cauchy-Riemann
systems. The first defines so-called $p$-analytic functions \cite{Polozhy}
(see also \cite{APFT}) and has the following form
\begin{equation}
\varphi_{x}=\frac{1}{p}\psi_{y},\qquad\varphi_{y}=-\frac{1}{p}\psi_{x}
\label{p-anal}%
\end{equation}
where $p$ is a given positive function of two real variables $x$ and $y$ which
is supposed to be continuously differentiable; $\varphi$ and $\psi$ are
real-valued continuously differentiable functions. If $\varphi$ and $\psi$ are
solutions of this system in $\Omega$, then the complex function $\omega
=\varphi+i\psi$ of a complex variable $z=x+iy$ is said to be $p$-analytic in
$\Omega$.

The second system considered here represents a special case of a general Vekua
equation (see, e.g., \cite{Vekua}) and sometimes is referred to as the main
Vekua equation \cite{APFT}. It has the form
\begin{equation}
W_{\overline{z}}=\frac{f_{\overline{z}}}{f}\overline{W}\text{\qquad in }%
\Omega\label{VekuamainA}%
\end{equation}
where the subindex $\overline{z}$ means the application of the operator
$\partial_{\overline{z}}:=\frac{1}{2}(\partial_{x}+i\partial_{y})$, $W$ is a
continuously differentiable complex valued function, $f$ is a positive
function of $x$ and $y$, twice continuously differentiable, which will be
supposed to be a particular solution of some stationary Schr\"{o}dinger
equation
\begin{equation}
\left(  -\Delta+q\right)  u=0\text{\qquad in }\Omega, \label{Schrod}%
\end{equation}
that is $q=\Delta f/f$.

Systems (\ref{p-anal}) and (\ref{VekuamainA}) are equivalent
\cite{KrTransplant}, \cite{APFT} in the following sense. Denote
\[
\mathcal{V}:=\partial_{\overline{z}}-\frac{f_{\overline{z}}}{f}C
\]
where $C$ is the operator of complex conjugation. We suppose that $p=f^{2}$
and introduce the operator%

\[
\Pi:=f\partial_{\overline{z}}P^{+}+\frac{1}{f}\partial_{\overline{z}}P^{-},
\]
where $P^{\pm}:=\frac{1}{2}(I\pm C)$ and $I$ is the identity operator. We have
that the equation
\begin{equation}
\Pi\omega=0 \label{PolozhOperatorForm}%
\end{equation}
is equivalent to the system%
\begin{equation}
\varphi_{x}=\frac{1}{f^{2}}\psi_{y},\qquad\varphi_{y}=-\frac{1}{f^{2}}\psi_{x}
\label{Polozhy4}%
\end{equation}
where $\varphi=\operatorname*{Re}\omega$ and $\psi=\operatorname*{Im}\omega$.

Denote
\[
B:=fP^{+}+\frac{1}{f}P^{-}.
\]
Then it is easy to see that
\[
B^{-1}=\frac{1}{f}P^{+}+fP^{-}.
\]

\begin{proposition}
\label{PrRelation}\cite{KrTransplant}%
\[
\mathcal{V}B=\Pi.
\]

\end{proposition}

\begin{remark}
From proposition \ref{PrRelation} we have also that%
\[
\mathcal{V}=\Pi B^{-1}.
\]

\end{remark}

Thus, application of the operator $B$ or $B^{-1}$ respectively allows us to
establish a direct relation between the results corresponding to
(\ref{p-anal}) and (\ref{VekuamainA}).

The following factorization of the Schr\"{o}dinger operator will be used.

\begin{theorem}
\cite{Krpseudoan} \label{Thfact}Let $f$ be a positive in $\Omega$ particular
solution of (\ref{Schrod}). Then for any real valued function $\varphi\in
C^{2}(\Omega)$ the following equalities hold%
\begin{equation}
\frac{1}{4}\left(  \Delta-\nu\right)  \varphi=\left(  \partial_{\overline{z}%
}+\frac{f_{z}}{f}C\right)  \left(  \partial_{z}-\frac{f_{z}}{f}C\right)
\varphi=\left(  \partial_{z}+\frac{f_{\overline{z}}}{f}C\right)  \left(
\partial_{\overline{z}}-\frac{f_{\overline{z}}}{f}C\right)  \varphi.
\label{fact}%
\end{equation}

\end{theorem}

An immediate corollary of this theorem is the fact that if $W$ is a solution
of (\ref{VekuamainA}) then its real part $W_{1}$ is necessarily a solution of
(\ref{Schrod}), meanwhile its imaginary part $W_{2}$ is a solution of the
following Schr\"{o}dinger equation%
\begin{equation}
-\Delta W_{2}+q_{1}W_{2}=0\qquad\text{in }\Omega\label{Schr2}%
\end{equation}
where $q_{1}=2(\nabla f)^{2}/f^{2}-q$ and $(\nabla f)^{2}=f_{x}^{2}+f_{y}^{2}$
(see \cite{Krpseudoan} and \cite{APFT}). Moreover, given $W_{1}$, the
corresponding $W_{2}$ can be easily constructed and vice versa. In order to
formulate this result we need to introduce the following notation. Note that
the operator $\partial_{\overline{z}}$ applied to a real-valued function
$\phi$ can be regarded as a kind of gradient, and if we know that
$\partial_{\overline{z}}\phi=\Phi$ in a whole complex plane or in a convex
domain, where $\Phi=\Phi_{1}+i\Phi_{2}$ is a given complex valued function
such that its real part $\Phi_{1}$ and imaginary part $\Phi_{2}$ satisfy the
equation
\begin{equation}
\partial_{y}\Phi_{1}-\partial_{x}\Phi_{2}=0, \label{casirot}%
\end{equation}
then we can reconstruct $\phi$ up to an arbitrary real constant $c$ in the
following way%
\begin{equation}
\phi(x,y)=2\left(  \int_{x_{0}}^{x}\Phi_{1}(\eta,y)d\eta+\int_{y_{0}}^{y}%
\Phi_{2}(x_{0},\xi)d\xi\right)  +c \label{Antigr}%
\end{equation}
where $(x_{0},y_{0})$ is an arbitrary fixed point in the domain of interest.
Note that this formula can be easily extended to any simply connected domain
by considering the integral along an arbitrary rectifiable curve $\Gamma$
leading from $(x_{0},y_{0})$ to $(x,y)$%
\[
\phi(x,y)=2\left(  \int_{\Gamma}\Phi_{1}dx+\Phi_{2}dy\right)  +c.
\]
By $\overline{A}$ we denote this integral operator:%
\[
\overline{A}[\Phi](x,y)=2\left(  \int_{x_{0}}^{x}\Phi_{1}(\eta,y)d\eta
+\int_{y_{0}}^{y}\Phi_{2}(x_{0},\xi)d\xi\right)  +c.
\]
Thus if $\Phi$ satisfies (\ref{casirot}), there exists a family of real valued
functions $\phi$ such that $\partial_{\overline{z}}\varphi=\Phi$, given by the
formula $\phi=\overline{A}[\Phi]$.

\bigskip

\begin{theorem}
\label{PrTransform}\cite{Krpseudoan} Let $W_{1}$ be a real valued solution of
(\ref{Schrod}) in a simply connected domain $\Omega$. Then the real valued
function $W_{2},$ solution of (\ref{Schr2}) such that $W=W_{1}+iW_{2}$ is a
solution of (\ref{VekuamainA}), is constructed according to the formula%
\begin{equation}
W_{2}=f^{-1}\overline{A}(if^{2}\partial_{\overline{z}}(f^{-1}W_{1})).
\label{transfDarboux}%
\end{equation}

Given a solution $W_{2}$ of (\ref{Schr2}), the corresponding solution $W_{1}$
of (\ref{Schrod}) such that $W=W_{1}+iW_{2}$ is a solution of
(\ref{VekuamainA}), is constructed as follows%
\begin{equation}
W_{1}=-f\overline{A}(if^{-2}\partial_{\overline{z}}(fW_{2})).
\label{transfDarbouxinv}%
\end{equation}

\end{theorem}

\bigskip

\begin{remark}
When in (\ref{Schrod}) $q\equiv0$ and $f\equiv1$, equalities
(\ref{transfDarboux}) and (\ref{transfDarbouxinv}) turn into the well known
formulas in complex analysis for constructing conjugate harmonic functions.
\end{remark}

We will need some definitions and results from Bers' pseudoanalytic function
theory \cite{Berskniga} concerning solutions of the general Vekua equation%
\begin{equation}
W_{\overline{z}}=a_{(F,G)}W+b_{(F,G)}\overline{W}, \label{VekuaGen}%
\end{equation}
where we will suppose that $a_{(F,G)}$ and $\ b_{(F,G)}$ are continuously
differentiable complex functions. A couple of solutions of (\ref{VekuaGen}) in
$\Omega$, $F$ and $G$ satisfying the inequality $\operatorname{Im}%
(\overline{F}G)>0$ form a so-called generating pair of the Vekua equation.
Every complex function $W$ defined in $\Omega$ admits the unique
representation $W=\phi F+\psi G$ where the functions $\phi$ and $\psi$ are
real valued. Sometimes it is convenient to associate with the function $W$ the
function $\omega=\phi+i\psi$. The correspondence between $W$ and $\omega$ is one-to-one.

The following expressions are known as characteristic coefficients of the pair
$(F,G)$
\[
a_{(F,G)}=-\frac{\overline{F}G_{\overline{z}}-F_{\overline{z}}\overline{G}%
}{F\overline{G}-\overline{F}G},\qquad b_{(F,G)}=\frac{FG_{\overline{z}%
}-F_{\overline{z}}G}{F\overline{G}-\overline{F}G},
\]

\[
A_{(F,G)}=-\frac{\overline{F}G_{z}-F_{z}\overline{G}}{F\overline{G}%
-\overline{F}G},\qquad B_{(F,G)}=\frac{FG_{z}-F_{z}G}{F\overline{G}%
-\overline{F}G}.
\]

For solutions of (\ref{VekuaGen}) the following operation is introduced,
called the $(F,G)$-derivative and denoted as $\overset{\cdot}{W}%
=\frac{d_{(F,G)}W}{dz}$:
\begin{equation}
\overset{\cdot}{W}=W_{z}-A_{(F,G)}W-B_{(F,G)}\overline{W}=\phi_{z}F+\psi_{z}G.
\label{FGder}%
\end{equation}

The inverse operation is introduced as follows.

\begin{definition}
\label{DefAdjoint}Let $(F,G)$ be a generating pair. Its adjoint generating
pair $(F,G)^{\ast}=(F^{\ast},G^{\ast})$ is defined by the formulas%
\[
F^{\ast}=-\frac{2\overline{F}}{F\overline{G}-\overline{F}G},\qquad G^{\ast
}=\frac{2\overline{G}}{F\overline{G}-\overline{F}G}.
\]

\end{definition}

The $(F,G)$-integral is defined as follows
\begin{equation}
\int_{\Gamma}Wd_{(F,G)}z=F(z_{1})\operatorname{Re}\int_{\Gamma}G^{\ast
}Wdz+G(z_{1})\operatorname{Re}\int_{\Gamma}F^{\ast}Wdz \label{pseudointegral}%
\end{equation}
where $\Gamma$ is a rectifiable curve leading from $z_{0}$ to $z_{1}$.

If $W=\phi F+\psi G$ is a solution of (\ref{VekuaGen}) where $\phi$ and $\psi$
are real valued functions then
\begin{equation}
\int_{z_{0}}^{z}\overset{\cdot}{W}d_{(F,G)}z=W(z)-\phi(z_{0})F(z)-\psi
(z_{0})G(z), \label{FGAnt}%
\end{equation}
and as $\overset{\cdot}{F}=\overset{}{\overset{\cdot}{G}=}0$, this integral is
path-independent and represents the $(F,G)$-antiderivative of $\overset{\cdot
}{W}$.

The $(F,G)$-derivative $\overset{\cdot}{W}$ is a solution of another Vekua
equation with some other coefficients $(a_{1},b_{1})$ and possessing another
generating pair $(F_{1},G_{1})$ called a successor of $(F,G)$.

\begin{definition}
\label{DefSeq}A sequence of generating pairs $\left\{  (F_{m},G_{m})\right\}
$, $m=0,\pm1,\pm2,\ldots$ , is called a generating sequence if $(F_{m+1}%
,G_{m+1})$ is a successor of $(F_{m},G_{m})$. If $(F_{0},G_{0})=(F,G)$, we say
that $(F,G)$ is embedded in $\left\{  (F_{m},G_{m})\right\}  $.
\end{definition}

Let $W$ be an $(F,G)$-pseudoanalytic function. Using a generating sequence in
which $(F,G)$ is embedded we can define the higher derivatives of $W$ by the
recursion formula%
\[
W^{[0]}=W;\qquad W^{[m+1]}=\frac{d_{(F_{m},G_{m})}W^{[m]}}{dz},\quad
m=1,2,\ldots\text{.}%
\]

The notion of a generating sequence leads to the concept of formal powers.

\begin{definition}
Each formal power $Z^{(n)}(a,z_{0};z)$ corresponding to the Vekua equation
(\ref{VekuaGen}), with some exponent $n\in\mathbb{Z}$, $a$ being a complex
number, $z_{0}$ a point in $\Omega$, is a solution of (\ref{VekuaGen}) in the
whole domain $\Omega$, such that
\[
\lim_{z\rightarrow z_{0}}\frac{Z^{(n)}(a,z_{0};z)}{a(z-z_{0})^{n}}=1.
\]
That is $Z^{(n)}(a,z_{0};z)$ is a solution of (\ref{VekuaGen}) possessing a
zero or a pole of order $n$ depending on the sign of $n$, and for $n=0$ it
takes the value $a$ at $z_{0}$. \label{defFormalPowers}
\end{definition}

The nonnegative formal powers ($n\geq0$) can be defined also in the following
recursive way.

\begin{definition}
\label{DefFormalPower}The formal power $Z_{m}^{(0)}(a,z_{0};z)$ with center at
$z_{0}\in\Omega$, coefficient $a$ and exponent $0$ is defined as the linear
combination of the generators $F_{m}$, $G_{m}$ with real constant coefficients
$\lambda$, $\mu$ chosen so that $\lambda F_{m}(z_{0})+\mu G_{m}(z_{0})=a$. The
formal powers with exponents $n=1,2,\ldots$ are defined by the recursion
formula%
\begin{equation}
Z_{m}^{(n)}(a,z_{0};z)=n\int_{z_{0}}^{z}Z_{m+1}^{(n-1)}(a,z_{0};\zeta
)d_{(F_{m},G_{m})}\zeta. \label{recformula}%
\end{equation}

\end{definition}

This definition implies the following properties.

\begin{enumerate}
\item $Z_{m}^{(n)}(a,z_{0};z)$ is an $(F_{m},G_{m})$-pseudoanalytic function
of $z$, that is, it is a solution of the Vekua equation $w_{\overline{z}%
}=a_{m}w+b_{m}\overline{w}$ possessing a generating pair $(F_{m},G_{m})$.

\item If $a^{\prime}$ and $a^{\prime\prime}$ are real constants, then
$Z_{m}^{(n)}(a^{\prime}+ia^{\prime\prime},z_{0};z)=a^{\prime}Z_{m}%
^{(n)}(1,z_{0};z)+a^{\prime\prime}Z_{m}^{(n)}(i,z_{0};z).$

\item The formal powers satisfy the differential relations%
\begin{equation}
\frac{d_{(F_{m},G_{m})}Z_{m}^{(n)}(a,z_{0};z)}{dz}=nZ_{m+1}^{(n-1)}%
(a,z_{0};z). \label{diffFormalPowers}%
\end{equation}

\item The asymptotic formulas
\[
Z_{m}^{(n)}(a,z_{0};z)\sim a(z-z_{0})^{n},\quad z\rightarrow z_{0}%
\]
hold.
\end{enumerate}

Moreover, the system of all formal powers $\left\{  Z_{0}^{(n)}(a,z_{0}%
;z)\right\}  _{n=0}^{\infty}$ represents a complete system of solutions of
(\ref{VekuaGen}) in the following sense. We will omit the subindex $0$ when a
formal power corresponds to $(F,G)$, that is $Z^{(n)}(a,z_{0};z):=Z_{0}%
^{(n)}(a,z_{0};z)$.

\begin{theorem}
\cite{BersFormalPowers}\label{ThRunge} A solution of (\ref{VekuaGen}) defined
in a bounded simply connected domain can be expanded into a normally
convergent series of formal polynomials (linear combinations of formal powers
with positive exponents).
\end{theorem}

Moreover, the following stronger result is valid.

\begin{theorem}
\cite{Menke}\label{ThMenke} Let $W$ be a solution of (\ref{VekuaGen}) in a
domain $\Omega$ bounded by a Jordan curve and satisfy the H\"{o}lder condition
on $\partial\Omega$ with the exponent $\alpha$ ($0<\alpha\leq1$). Then for any
$\varepsilon>0$ and any natural $n$ there exists a pseudopolynomial of order
$n$ satisfying the inequality
\[
\left\vert W(z)-P_{n}(z)\right\vert \leq\frac{\operatorname*{Const}}%
{n^{\alpha-\varepsilon}}\qquad\text{for any }z\in\overline{\Omega}%
\]
where the constant does not depend on $n$, but only on $\varepsilon$.
\end{theorem}

With the aid of these results concerning pseudoanalytic formal powers and of
the relation between solutions of the main Vekua equation to the
Schr\"{o}dinger equation corresponding completeness results were obtained for
solutions of the Schr\"{o}dinger equation as, e.g., the following statement.

\begin{theorem}
\cite{KrJPhys06}\label{ThRungeSchr} An arbitrary solution of (\ref{Schrod})
defined in a bounded simply connected domain $\Omega$ where there exists a
positive particular solution $f\in C^{1}(\overline{\Omega})$ of (\ref{Schrod})
can be expanded into a normally convergent series of real parts of formal polynomials.
\end{theorem}

As was mentioned before besides positive formal powers also the negative were
defined by L. Bers (see \cite{Berskniga}). First of all, the existence of the
generalized Cauchy kernel was proved, that is the existence of a solution $w$
of (\ref{VekuaGen}) in $\Omega\setminus\left\{  z_{0}\right\}  $ which
satisfies the relation%
\begin{equation}
\lim_{z\rightarrow z_{0}}\frac{w(z)}{a(z-z_{0})^{-1}}=1
\label{condCauchykernel}%
\end{equation}
where $a$ is any complex number. This function is denoted as follows%
\[
w(z)=Z^{(-1)}(a,z_{0},z).
\]
The negative formal powers $Z^{(-n)}$ for\ $n=2,3,\ldots$, are constructed
using the recursive differential relations like (\ref{diffFormalPowers}).

With the aid of the positive and negative formal powers a whole theory of
pseudoanalytic functions was developed including Taylor and Laurent series,
and their numerous properties similar to the properties of their special cases
corresponding to the usual analytic functions. The generalized Cauchy kernel
$Z^{(-1)}(\alpha,z_{0},z)$ makes it possible to prove a generalization of the
Cauchy integral formula \cite{Berskniga}, see also \cite{APFT}.

Thus, an important problem is to find the way to construct the formal powers
explicitly. This is the main subject of this paper.

\section{The transplant operator\label{SectTransplant}}

In this section we define and study the main tool of this paper called the
transplant operator. It was introduced in \cite{KrTransplant} and used for
constructing Cauchy kernels and Cauchy integral representations for an
important subclass of $p$-analytic functions,- the $x^{k}$-analytic functions.
Let us describe the main idea behind this concept.

Let both $f$ and $g$ be positive solutions of (\ref{Schrod}) in $\Omega$.
Together with the main Vekua equation (\ref{VekuamainA}) we consider the main
Vekua equation corresponding to $g$:%
\begin{equation}
w_{\overline{z}}=\frac{g_{\overline{z}}}{g}\overline{w}\text{\qquad in }%
\Omega. \label{VekuamainForg}%
\end{equation}
We have that both $\operatorname*{Re}W$ (where $W$ is a solution of
(\ref{VekuamainA})) and $\operatorname*{Re}w$ satisfy (\ref{Schrod}) in
$\Omega$, meanwhile $\operatorname*{Im}W$ and $\operatorname*{Im}w$ satisfy in
general different Schr\"{o}dinger equations%
\begin{equation}
(-\Delta+q_{1})\operatorname*{Im}W=0\text{\qquad in }\Omega\label{schroq1}%
\end{equation}
and
\begin{equation}
(-\Delta+q_{2})\operatorname*{Im}w=0\text{\qquad in }\Omega\label{schroq2}%
\end{equation}
where $q_{1}=2(\nabla f)^{2}/f^{2}-q$ and $q_{2}=2(\nabla g)^{2}/g^{2}-q$.

Now we introduce an operator which transforms solutions of (\ref{VekuamainA})
into solutions of (\ref{VekuamainForg}) acting in the following way%
\begin{equation}
T_{f,g}[W]=P^{+}W+ig^{-1}\overline{A}[ig^{2}\partial_{\overline{z}}%
(g^{-1}P^{+}W)]. \label{transplantop}%
\end{equation}
Its application makes the imaginary part of a solution of (\ref{VekuamainA})
drop out and be substituted by an imaginary part constructed according to
theorem \ref{PrTransform} in such a way that after this \textquotedblleft
transplant\textquotedblright\ operation the new complex function
$w=T_{f,g}[W]$ becomes a solution of (\ref{VekuamainForg}). This is why we
call the operator $T_{f,g}$ the transplant operator.

Assigning a fixed value in a certain point of the domain of interest to the
result of application of $\overline{A}$ we obtain an invertible one-to-one map
establishing a relation between solutions of (\ref{VekuamainA}) and
(\ref{VekuamainForg}). The inverse to $T_{f,g}$ is given by the expression%
\[
T_{f,g}^{-1}[w]=T_{g,f}[w]=P^{+}w+if^{-1}\overline{A}[if^{2}\partial
_{\overline{z}}(f^{-1}P^{+}w)].
\]

Let us denote the formal powers corresponding to (\ref{VekuamainA}) and
(\ref{VekuamainForg}) by $Z_{f}^{(n)}(a,z_{0},z)$ and $Z_{g}^{(n)}(a,z_{0},z)$
respectively. In the following we establish a useful property of the
transplant operator. Namely, that it allows one to transform an $n$-th formal
power to an $n$-th formal power. We will consider the case of positive and
negative formal powers separately.

Let $W$ be a solution of (\ref{VekuamainA}) such that
\begin{equation}
\lim_{z\rightarrow z_{0}}\frac{W(z)}{a(z-z_{0})^{n}}=1 \label{asympt}%
\end{equation}
for some $z_{0}\in\Omega$, $n\in\mathbb{N}$ and a complex number $a$. That is
$W$ is a formal power $Z_{f}^{(n)}(a,z_{0},z)$ corresponding to
(\ref{VekuamainA}). As before, we denote $W_{1}=\operatorname*{Re}W$ and
$W_{2}=\operatorname{Im}W$ and due to theorem \ref{PrTransform} we have the
equality (\ref{transfDarboux}). As $W$ has a zero at $z_{0}$ it is convenient
to write $\overline{A}[\Phi]$ where $\Phi=if^{2}\partial_{\overline{z}}%
(f^{-1}W_{1})$ as follows%

\begin{equation}
\overline{A}[\Phi](z)=2\left(  \int_{\Gamma}\Phi_{1}dx+\Phi_{2}dy\right)
\label{Abar}%
\end{equation}
where $\Gamma$ is a rectifiable curve leading from $z_{0}$ to $z$. That is we
fix $z_{0}$ as an initial point for integration in $\overline{A}$.

Now consider
\begin{equation}
\omega_{2}=g^{-1}\overline{A}(ig^{2}\partial_{\overline{z}}(g^{-1}W_{1}))
\label{omega2}%
\end{equation}
where again $z_{0}$ is used as an initial point for integration. We are
interested in the limit%
\begin{equation}
\lim_{z\rightarrow z_{0}}\frac{W_{2}(z)}{\omega_{2}(z)}=c\lim_{z\rightarrow
z_{0}}\frac{\psi_{f}(z)}{\psi_{g}(z)} \label{limit}%
\end{equation}
where $c:=f^{-1}(z_{0})/g^{-1}(z_{0})$,
\begin{equation}
\psi_{f}:=\overline{A}(if^{2}\partial_{\overline{z}}(f^{-1}W_{1}))
\label{psif}%
\end{equation}
and
\begin{equation}
\psi_{g}:=\overline{A}(ig^{2}\partial_{\overline{z}}(g^{-1}W_{1})).
\label{psig}%
\end{equation}
In order to prove its existence and evaluate it let us consider any direction
in the plane defined by a vector $\mathbf{d}=(d_{1},d_{2})^{T}$, and assume
that $z$ tends to $z_{0}$ along the corresponding path, that is we consider
the following limit%
\[
\lim_{t\rightarrow0}\frac{\psi_{f}(x_{0}+td_{1},y_{0}+td_{2})}{\psi_{g}%
(x_{0}+td_{1},y_{0}+td_{2})}.
\]
By definition, $\psi_{f}(z_{0})=\psi_{g}(z_{0})=0$, and hence to evaluate this
limit we can make use of the l'Hospital rule which here gives us%
\[
\lim_{t\rightarrow0}\frac{\psi_{f}(x_{0}+td_{1},y_{0}+td_{2})}{\psi_{g}%
(x_{0}+td_{1},y_{0}+td_{2})}=\frac{\frac{\partial\psi_{f}(x_{0},y_{0}%
)}{\partial\mathbf{d}}}{\frac{\partial\psi_{g}(x_{0},y_{0})}{\partial
\mathbf{d}}}=\frac{\left\langle \nabla\psi_{f}(x_{0},y_{0}),\mathbf{d}%
\right\rangle }{\left\langle \nabla\psi_{g}(x_{0},y_{0}),\mathbf{d}%
\right\rangle }%
\]
where $\left\langle \cdot,\cdot\right\rangle $ denotes the usual scalar
product of two vectors. Let us note that the last expression can be written in
a complex-analytic form as follows%
\[
\frac{\operatorname{Re}\left(  \partial_{\overline{z}}\psi_{f}(z_{0}%
)\cdot(d_{1}-id_{2})\right)  }{\operatorname{Re}\left(  \partial_{\overline
{z}}\psi_{g}(z_{0})\cdot(d_{1}-id_{2})\right)  }.
\]
We recall that $\psi_{f}$ and $\psi_{g}$ are defined by (\ref{psif}) and
(\ref{psig}) respectively. Thus we have
\begin{align*}
\lim_{t\rightarrow0}\frac{\psi_{f}(x_{0}+td_{1},y_{0}+td_{2})}{\psi_{g}%
(x_{0}+td_{1},y_{0}+td_{2})}  &  =\frac{\operatorname{Re}\left(
if^{2}\partial_{\overline{z}}(f^{-1}W_{1})\cdot(d_{1}-id_{2})\right)
}{\operatorname{Re}\left(  ig^{2}\partial_{\overline{z}}(g^{-1}W_{1}%
)\cdot(d_{1}-id_{2})\right)  }\\
&  =\frac{1}{c^{2}}\frac{\operatorname{Re}\left(  \partial_{\overline{z}%
}(f^{-1}W_{1})\cdot(d_{2}+id_{1})\right)  }{\operatorname{Re}\left(
\partial_{\overline{z}}(g^{-1}W_{1})\cdot(d_{2}+id_{1})\right)  }\\
&  =\frac{1}{c^{2}}\frac{\operatorname{Re}\left(  (\partial_{\overline{z}%
}f^{-1}(z_{0})W_{1}(z_{0})+f^{-1}(z_{0})\partial_{\overline{z}}W_{1}%
(z_{0}))\cdot(d_{2}+id_{1})\right)  }{\operatorname{Re}\left(  (\partial
_{\overline{z}}g^{-1}(z_{0})W_{1}(z_{0})+g^{-1}(z_{0})\partial_{\overline{z}%
}W_{1}(z_{0}))\cdot(d_{2}+id_{1})\right)  }.
\end{align*}
Now we use the fact that $W_{1}(z_{0})=0$ as well as once more that $f$ and
$g$ are positive and obtain that%
\[
\lim_{t\rightarrow0}\frac{\psi_{f}(x_{0}+td_{1},y_{0}+td_{2})}{\psi_{g}%
(x_{0}+td_{1},y_{0}+td_{2})}=\frac{1}{c}%
\]
for any direction $\mathbf{d}$. Thus, the limit (\ref{limit}) exists and
$\lim_{z\rightarrow z_{0}}\frac{W_{2}(z)}{\omega_{2}(z)}=1$. Consequently we
obtain that the function $W_{1}+i\omega_{2}$ satisfies the asymptotic relation
(\ref{asympt}) as well and represents a formal power $Z_{g}^{(n)}(a,z_{0},z)$
corresponding to (\ref{VekuamainForg}).

Now let us consider negative formal powers. We suppose that $W$ is a solution
of (\ref{VekuamainA}) such that
\begin{equation}
\lim_{z\rightarrow z_{0}}W(z)(z-z_{0})^{n}=a \label{asymptnegat}%
\end{equation}
for some $z_{0}\in\Omega$, $n\in\mathbb{N}$ and a complex number $a$.

As before, we denote $W_{1}=\operatorname*{Re}W$ and $W_{2}=\operatorname{Im}%
W$ and due to theorem \ref{PrTransform} we have the equality
(\ref{transfDarboux}) at any point $z\in\Omega$ distinct from $z_{0}$ where
$W$ has a pole of order $n$. The integration involved in $\overline{A}[\Phi]$
where $\Phi=if^{2}\partial_{\overline{z}}(f^{-1}W_{1})$ is done along any
rectifiable curve $\Gamma$ belonging to $\Omega,$ leading from $z_{1}$ to $z$
and not passing through $z_{0}$. Again we consider the function $\omega_{2}$
defined by (\ref{omega2}) where the integration is done in the same way as was
just explained. We are interested in the limit (\ref{limit}) and for this we
again consider any direction $\mathbf{d}$ and use the l%
\'{}%
Hospital rule as both functions tend to infinity at $z_{0}$:
\[
\lim_{t\rightarrow0}\frac{\psi_{f}(x_{0}+td_{1},y_{0}+td_{2})}{\psi_{g}%
(x_{0}+td_{1},y_{0}+td_{2})}=\frac{1}{c^{2}}\frac{\operatorname{Re}\left(
\partial_{\overline{z}}(f^{-1}W_{1})\cdot(d_{2}+id_{1})\right)  }%
{\operatorname{Re}\left(  \partial_{\overline{z}}(g^{-1}W_{1})\cdot
(d_{2}+id_{1})\right)  }.
\]
Here the reasoning we used before, in the case of positive formal powers, is
not already applicable. Nevertheless we note that the l%
\'{}%
Hospital rule can be applied to the obtained quotient in the opposite
direction. Namely, we have%
\begin{align*}
&  \frac{1}{c^{2}}\frac{\operatorname{Re}\left(  \partial_{\overline{z}%
}(f^{-1}W_{1})\cdot(d_{2}+id_{1})\right)  }{\operatorname{Re}\left(
\partial_{\overline{z}}(g^{-1}W_{1})\cdot(d_{2}+id_{1})\right)  }\\
&  =\frac{1}{c^{2}}\lim_{t\rightarrow0}\frac{f^{-1}\left(  x_{0}+td_{2}%
,y_{0}-td_{1}\right)  W_{1}\left(  x_{0}+td_{2},y_{0}-td_{1}\right)  }%
{g^{-1}\left(  x_{0}+td_{2},y_{0}-td_{1}\right)  W_{1}\left(  x_{0}%
+td_{2},y_{0}-td_{1}\right)  }=\frac{1}{c}.
\end{align*}

Thus we proved that with the aid of the transplant operator both positive and
negative formal powers corresponding to (\ref{VekuamainA}) and
(\ref{VekuamainForg}) can be transformed to each other. We formulate these
statements as the following theorems.

\begin{theorem}
Let $f$ and $g$ be real valued nonvanishing solutions of (\ref{Schrod}) in a
simply connected domain $\Omega\subset\mathbb{R}^{2}$. Let $z_{0}\in\Omega$,
$a\in\mathbb{C}$ and $Z_{f}^{(n)}(a,z_{0},z)$, $n\in\mathbb{N}$ be a formal
power associated with equation (\ref{VekuamainA}). Then the function
$Z_{g}^{(n)}(a,z_{0},z):=T_{f,g}\left[  Z_{f}^{(n)}(a,z_{0},z)\right]  $ is a
formal power of order $n$, with center at $z_{0}$ and coefficient $a$,
associated with equation (\ref{VekuamainForg}). Here $T_{f,g}$ is defined by
(\ref{transplantop}) with $\overline{A}$ being defined by (\ref{Abar}) where
as an initial point of integration is chosen $z_{0}$.
\end{theorem}

\begin{theorem}
Let $f$ and $g$ be real valued nonvanishing solutions of (\ref{Schrod}) in a
simply connected domain $\Omega\subset\mathbb{R}^{2}$. Let $z_{0}\in\Omega$,
$a\in\mathbb{C}$ and $Z_{f}^{(-n)}(a,z_{0},z)$, $n\in\mathbb{N}$ be a formal
power associated with equation (\ref{VekuamainA}). Then the function
$Z_{g}^{(-n)}(a,z_{0},z):=T_{f,g}\left[  Z_{f}^{(-n)}(a,z_{0},z)\right]  $ is
a formal power of order $-n$, with center at $z_{0}$ and coefficient $a$,
associated with equation (\ref{VekuamainForg}). Here $T_{f,g}$ is defined by
(\ref{transplantop}) with $\overline{A}$ being defined by (\ref{Abar}) where
$\Gamma$ is any rectifiable curve belonging to $\Omega$, leading from $z_{1}$
to $z$ and not passing through $z_{0}$.
\end{theorem}

\section{Construction of positive formal
powers\label{SectPositiveFormalPowers}}

As we have shown in the previous section the transplant operator allows us to
transform positive and negative formal powers of one main Vekua equation, say
(\ref{VekuamainA}), into formal powers of the same order of another main Vekua
equation, say (\ref{VekuamainForg}), when the coefficients $f$ and $g$ are
solutions of the same Schr\"{o}dinger equation (\ref{Schrod}). This
observation leads to a substantial extension of the class of Vekua equations
and of systems of the form (\ref{p-anal}) for which a generating sequence and
a complete system of formal powers can be obtained. Suppose we are interested
in solving a Vekua equation of the form (\ref{VekuamainForg}) or a system
describing $p$-analytic functions (\ref{p-anal}) with $p=g^{2}$. Then the
first step is to look for a \textquotedblleft simplest\textquotedblright%
\ solution $f$ of the equation (\ref{Schrod}) where $q=\Delta g/g$, such that
for the corresponding main Vekua equation (\ref{VekuamainA}) a generating
sequence and hence a system of formal powers can be constructed. Then
application of the transplant operator gives a system of formal powers for
(\ref{VekuamainForg}) and (\ref{p-anal}) as well as a corresponding generating sequence.

As an example, let us consider two positive solutions $f=y^{2}$ and
$g=\displaystyle\frac{1+xy^{3}}{y}$ of the Schr\"{o}dinger equation
(\ref{Schrod}) with potential $q=2/y^{2}$ in the domain $\Omega
=\{(x,y)\ |\ x>0\ \mbox{and}\ y>0\}$. In this case, the Vekua equations
(\ref{VekuamainA}) and (\ref{VekuamainForg}) take, respectively, the form
\begin{equation}
W_{\overline{z}}=\frac{i}{y}\overline{W}\text{\qquad in }\Omega,
\label{MainVekuaf}%
\end{equation}
and
\begin{equation}
w_{\overline{z}}=\frac{y^{4}+i(2xy^{3}-1)}{2y(1+xy^{3})}\overline
{w}\text{\qquad in }\Omega. \label{MainVekuag}%
\end{equation}

For the calculations given below we used Maple. Let us first calculate the
formal powers $Z_{f}^{(n)}(a,z_{0};z)$ of orders $n=0,1,2$ for the Vekua
equation (\ref{MainVekuaf}) with generating pair $(F,G)=(f,i/f)$ where
$z_{0}:=x_{0}+iy_{0}$ and $(x_{0},y_{0})\in\Omega$. Using property~2 following
definition~\ref{DefFormalPower}, we are considering $Z_{f}^{(n)}(1,z_{0};z)$
and $Z_{f}^{(n)}(i,z_{0};z)$. By definition \ref{DefFormalPower} we have
$Z_{f}^{(0)}(1,z_{0};z)=\lambda F(z)+\mu G(z)$ and $Z_{f}^{(0)}(i,z_{0}%
;z)=\lambda^{\prime}F(z)+\mu^{\prime}G(z)$ where the constants $(\lambda
,\mu),(\lambda^{\prime},\mu^{\prime})$ are defined by $\lambda F(z_{0})+\mu
G(z_{0})=1$ and $\lambda^{\prime}F(z_{0})+\mu^{\prime}G(z_{0})=i$. We find
$(\lambda,\mu)=(1/y_{0}^{2},0)$ and $(\lambda^{\prime},\mu^{\prime}%
)=(0,y_{0}^{2})$ such that
\[
Z_{f}^{(0)}(1,z_{0};z)=\left(  \frac{y}{y_{0}}\right)  ^{2}\quad
\quad\mbox{ and }\quad\quad Z_{f}^{(0)}(i,z_{0};z)=i\left(  \frac{y_{0}}%
{y}\right)  ^{2}.
\]
In order to construct $Z_{f}^{(1)}(\alpha,z_{0};z)$ for $\alpha=1,i$ from
formula (\ref{recformula}) we need first $Z_{f,1}^{(0)}(\alpha,z_{0};z)$.
However, for $f$ depending only on $y$ it is shown (see \cite{Berskniga,
APFT}) that $(F_{m},G_{m})=(F,G)$ for $m=0,\pm1,\pm2,\ldots$. Therefore we
have
\[
Z_{f,m}^{(n)}(\alpha,z_{0};z)=Z_{f}^{(n)}(\alpha,z_{0};z)\quad\quad
\mbox{ for }\quad\quad\alpha=1,i\ \ \mbox{ and }\ \ m=0,\pm1,\pm2,\ldots
\]
so that formula (\ref{recformula}) gives us
\[
Z_{f}^{(1)}(\alpha,z_{0};z)=\int_{z_{0}}^{z}Z_{f}^{(0)}(\alpha,z_{0}%
;\zeta)d_{(F,G)}\zeta,\quad\quad\alpha=1,i.
\]
We calculate these two integrals using (\ref{pseudointegral}) where $F^{\ast
}=-if$ and $G^{\ast}=1/f$. Defining $\zeta:=\xi+i\eta$, we obtain
\begin{align*}
Z_{f}^{(1)}(1,z_{0};z)  &  =y^{2}\mbox{Re}\displaystyle\int_{z_{0}}%
^{z}\displaystyle\frac{d\zeta}{y_{0}^{2}}-\displaystyle\frac{i}{y^{2}%
}\mbox{Re}\displaystyle\int_{z_{0}}^{z}\displaystyle\frac{i\eta^{4}}{y_{0}%
^{2}}d\zeta\\
&  =(x-x_{0})\left(  \displaystyle\frac{y}{y_{0}}\right)  ^{2}%
+\displaystyle\frac{i}{5}\frac{y^{5}-y_{0}^{5}}{(y_{0}y)^{2}}%
\end{align*}
and
\begin{align*}
Z_{f}^{(1)}(i,z_{0};z)  &  =y^{2}\mbox{Re}\displaystyle\int_{z_{0}}%
^{z}\displaystyle\frac{iy_{0}^{2}}{\eta^{4}}d\zeta+\displaystyle\frac{i}%
{y^{2}}\mbox{Re}\displaystyle\int_{z_{0}}^{z}y_{0}^{2}d\zeta\\
&  =-\displaystyle\frac{1}{3}\displaystyle\frac{y^{3}-y_{0}^{3}}{y_{0}%
y}+i(x-x_{0})\left(  \displaystyle\frac{y_{0}}{y}\right)  ^{2}.
\end{align*}

In a similar way we construct $Z_{f}^{(2)}(\alpha,z_{0};z)$ for $\alpha=1,i$
where we first need $Z_{f,1}^{(1)}(\alpha,z_{0};z)=Z_{f}^{(1)}(\alpha
,z_{0};z)$. From formula (\ref{recformula}) we obtain
\begin{align*}
Z_{f}^{(2)}(1,z_{0};z)  &  =2\displaystyle\int_{z_{0}}^{z}Z^{(1)}%
(1,z_{0};z)d_{(F,G)}\zeta\\
&  =2y^{2}\mbox{Re}\displaystyle\int_{z_{0}}^{z}\left[  (\xi-x_{0})\left(
\displaystyle\frac{\eta}{y_{0}}\right)  ^{2}+\frac{i}{5}\frac{\eta^{5}%
-y_{0}^{5}}{(y_{0}\eta)^{2}}\right]  \frac{d\zeta}{\eta^{2}}\\
&  -\displaystyle\frac{2i}{y^{2}}\mbox{Re}\displaystyle\int_{z_{0}}^{z}\left[
(\xi-x_{0})\left(  \displaystyle\frac{\eta}{y_{0}}\right)  ^{2}+\frac{i}%
{5}\frac{\eta^{5}-y_{0}^{5}}{(y_{0}\eta)^{2}}\right]  (i\eta^{2})d\zeta\\
&  =\displaystyle\frac{1}{15(y_{0}y)^{2}}\Big[\big(15(x-x_{0})^{2}y^{4}%
-3y^{6}+5y_{0}^{2}y^{4}-2y_{0}y\big)\\
&  +6i(x-x_{0})(y^{5}-y_{0}^{5})\Big]
\end{align*}
and
\begin{align*}
Z_{f}^{(2)}(i,z_{0};z)  &  =2\displaystyle\int_{z_{0}}^{z}Z^{(1)}%
(i,z_{0};z)d_{(F,G)}\zeta\\
&  =2y^{2}\mbox{Re}\displaystyle\int_{z_{0}}^{z}\left[  \frac{1}{3}\frac
{y_{0}^{3}-\eta^{3}}{y_{0}\eta}+i(\xi-x_{0})\left(  \displaystyle\frac{y_{0}%
}{\eta}\right)  ^{2}\right]  \frac{d\zeta}{\eta^{2}}\\
&  -\displaystyle\frac{2i}{y^{2}}\mbox{Re}\displaystyle\int_{z_{0}}^{z}\left[
\frac{1}{3}\frac{y_{0}^{3}-\eta^{3}}{y_{0}\eta}+i(\xi-x_{0})\left(
\displaystyle\frac{y_{0}}{\eta}\right)  ^{2}\right]  (i\eta^{2})d\zeta\\
&  =\displaystyle\frac{1}{15y_{0}y^{2}}\Big[-10(x-x_{0})y(y^{3}-y_{0}^{3})\\
&  +i\big(15y_{0}^{3}(x-x_{0})^{2}+5y_{0}^{3}y^{2}-2y^{5}-3y_{0}%
^{5}\big)\Big].
\end{align*}
We verify easily that $Z_{f}^{(n)}$ are indeed solutions of the Vekua equation
(\ref{MainVekuaf}). Moreover, $\mathrm{Re}\,Z_{f}^{(n)}$ are solutions of the
Schr\"{o}dinger equation (\ref{Schrod}) with $q=2/y^{2}$ and $\mathrm{Im}%
\,Z_{f}^{(n)}$ are solutions of the Schr\"{o}dinger equation (\ref{schroq1})
with $q_{1}=6/y^{2}$. Finally, we have (see definition \ref{defFormalPowers})
\[
\displaystyle\lim_{z\rightarrow z_{0}}\displaystyle\frac{Z_{f}^{(n)}%
(\alpha,z_{0};z)}{(z-z_{o})^{n}}=\alpha,\quad\quad\quad\alpha=1,i.
\]

Now in order to obtain the formal powers $Z_{g}^{(n)}(\alpha,z_{0}%
;z)=T_{f,g}\big[Z_{f}^{(n)}(\alpha,z_{0};z)\big]$ of Vekua equation
(\ref{MainVekuag}) let us apply the transplant operator to the constructed
formal powers $Z_{f}^{(n)}(\alpha,z_{0};z)$ of (\ref{MainVekuaf}) for $n=1,2$.

Since the formal powers $Z_{g}^{(0)}(\alpha,z_{0};z)$ can be easily calculated
using definition~\ref{DefFormalPower}, we are not using the transplant
operator in the particular case of formal powers of zero order. Hence, for the
generating pair $(F,G)=(g,i/g)$ we find
\[
Z_{g}^{(0)}(1,z_{0};z)=k_{0}\frac{1+xy^{3}}{y} \quad\mbox{and} \quad
Z_{g}^{(0)}(i,z_{0};z)=\frac{i}{k_{0}}\frac{y}{1+xy^{3}},
\]
where $k_{0}:=y_{0}/(1+x_{0}y_{0}^{3})$.

Now considering application of the transplant operator to $Z_{f}^{(n)}%
(\alpha,z_{0};z)$ for $n=1,2$ we obtain
\begin{align*}
Z_{g}^{(1)}(1,z_{0};z)  &  :=T_{f,g}\big[Z_{f}^{(1)}(1,z_{0};z)\big]\\
&  =\mathrm{Re}\,Z_{f}^{(1)}(1,z_{0};z)+ig^{-1}\overline{A}%
\Big[\displaystyle\frac{-3(x-x_{0})+i(y+x_{0}y^{4})}{2y_{0}^{2}}\Big],
\end{align*}
where
\begin{align*}
\overline{A}\Big[\displaystyle\frac{-3(x-x_{0})+i(y+x_{0}y^{4})}{2y_{0}^{2}%
}\Big]  &  =\displaystyle\int_{x_{0}}^{x}\displaystyle\frac{-3(\eta-x_{0}%
)}{{y_{0}^{2}}}d\eta+\displaystyle\int_{y_{0}}^{y}\frac{y+x_{0}y^{4}}%
{y_{0}^{2}}d\xi+c\\
&  =\displaystyle\frac{5(y^{2}-y_{0}^{2})+2x_{0}(y^{5}-y_{0}^{5}%
)-15(x-x_{0})^{2}}{10y_{0}^{2}}+c.
\end{align*}
Therefore, we have
\[
Z_{g}^{(1)}(1,z_{0};z)=(x-x_{0})\left(  \displaystyle\frac{y}{y_{0}}\right)
^{2}+i\displaystyle\frac{\big[5(y^{2}-y_{0}^{2})+2x_{0}(y^{5}-y_{0}%
^{5})-15(x-x_{0})^{2}\big]y}{10y_{0}^{2}(1+xy^{3})},
\]
where the arbitrary real constant $c$ was chosen equal to zero.

Similar calculations give us:
\begin{align*}
&  Z_{g}^{(1)}(i,z_{0};z)=-\displaystyle\frac{1}{3}\displaystyle\frac
{y^{3}-y_{0}^{3}}{y_{0}y}+i\frac{\big[30(x-x_{0})+15y_{0}^{3}(x^{2}-x_{0}%
^{2})-5y_{0}^{3}y^{2}+2y^{5}+3y_{0}^{5}\big]y}{30y_{0}(1+xy^{3})}\\
&  Z_{g}^{(2)}(1,z_{0};z)=\displaystyle\displaystyle\frac{15(x-x_{0})^{2}%
y^{4}-3y^{6}+5y_{0}^{2}y^{4}-2y_{0}^{5}y}{15(y_{0}y)^{2}}+\frac{iy}{105
y_{0}^{2}(1+xy^{3})}\Big[315x_{0}x(x-x_{0})\\
&  +105(x-x_{0})(y^{2}-y_{0}^{2})-105(x^{3}-x_{0}^{3})+21(x^{2}-x_{0}%
^{2})(y^{5}-y_{0}^{5})-7(y_{0}y)^{2}(y^{3}-y_{0}^{3})+3(y^{7}-y_{0}%
^{7})\Big]\\
&  Z_{g}^{(2)}(i,z_{0};z)=-\displaystyle\frac{2(x-x_{0})y(y^{3}-y_{0}^{3}%
)}{3y_{0}y^{2}}+\frac{i}{15y_{0}(1+xy^{3})}\Big[15(x-x_{0})^{2}y+10y_{0}%
^{3}x^{3}y\\
&  -15x_{0}y_{0}^{3}x^{2}y+5x_{0}^{3}y_{0}^{3}y-2x_{0}y^{6}-5y^{3}+5x_{0}%
y_{0}^{3}y^{3}-10y_{0}^{3}-3x_{0}y_{0}^{5}y+15y_{0}^{2}y\Big]
\end{align*}
These formal powers $Z_{g}^{(n)}$ are solutions of the Vekua equation
(\ref{MainVekuag}). Moreover, we also have that $\mathrm{Re}\,Z_{g}^{(n)}$ are
solutions of the Schr\"{o}dinger equation (\ref{Schrod}) with $q=2/y^{2}$ and
$\mathrm{Im}\,Z_{g}^{(n)}$ are solutions of the Schr\"{o}dinger equation
(\ref{schroq2}) with potential
\[
q_{2}=\frac{y(y^{5}+3x^{2}y^{3}-6x)}{(1+xy^{3})^{2}}.
\]
Finally, we can verify that $Z_{g}^{(n)}=T_{f,g}\big[Z_{f}^{(n)}\big]$ satisfy
the asymptotics of the formal powers when $z\rightarrow z_{0}$, i.e.
\[
\displaystyle\lim_{z\rightarrow z_{0}}\displaystyle\frac{Z_{g}^{(n)}%
(\alpha,z_{0};z)}{(z-z_{o})^{n}}=\alpha,\quad\quad\quad\alpha=1,i.
\]

\section{Construction of a generating sequence}

Meanwhile in the example considered in the previous section the generating
sequence for the equation (\ref{MainVekuaf}) is very simple and consists of
one generating pair only $(y^{2},i/y^{2})$, the generating sequence for the
related equation (\ref{MainVekuag}) is more complicated. However the procedure
based on the application of the transplant operator allows us to obtain a
generating sequence for a \textquotedblleft more complicated\textquotedblright%
\ main Vekua equation from a generating sequence corresponding to a
\textquotedblleft simpler\textquotedblright\ main Vekua equation. Here the
algorithm is following. First, using a generating sequence for equation
(\ref{VekuamainA}), which is assumed to be known, one can construct the
complete system of positive formal powers corresponding to (\ref{VekuamainA}).
Next, as was explained in the preceding two sections, application of the
transplant operator gives a complete system of positive formal powers for
equation (\ref{VekuamainForg}) where $g$ is related to $f$ via the
Schr\"{o}dinger equation (\ref{Schrod}). Finally, to obtain a generating
sequence for (\ref{VekuamainForg}) one can use property 3 of formal powers. We
illustrate this by the following scheme.

\noindent%
\begin{tabular}
[c]{cccccccc}%
${\tiny \vdots}$ &  &  &  &  &  &  & \\
$_{{\tiny (Z}^{(3)}{\tiny (1,z}_{0}{\tiny ;z),Z}^{(3)}{\tiny (i,z}%
_{0}{\tiny ;z))}}$ &  &  &  &  &  &  & \\
& $_{\overset{\frac{d_{(g,i/g)}}{dz}}{\searrow}}$ & ${\tiny \vdots}$ &  &  &
&  & \\
$_{{\tiny (Z}^{(2)}{\tiny (1,z}_{0}{\tiny ;z),Z}^{(2)}{\tiny (i,z}%
_{0}{\tiny ;z))}}$ &  & $_{{\tiny (Z}_{1}^{(2)}{\tiny (1,z}_{0}{\tiny ;z),Z}%
_{1}^{(2)}{\tiny (i,z}_{0}{\tiny ;z))}}$ &  &  &  &  & \\
& $_{\overset{\frac{d_{(g,i/g)}}{dz}}{\searrow}}$ &  & $_{_{\overset
{\frac{d_{(F_{1},G_{1})}}{dz}}{\searrow}}}$ & ${\tiny \vdots}$ &  &  & \\
$_{{\tiny (Z}^{(1)}{\tiny (1,z}_{0}{\tiny ;z),Z}^{(1)}{\tiny (i,z}%
_{0}{\tiny ;z))}}$ &  & $_{{\tiny (Z}_{1}^{(1)}{\tiny (1,z}_{0}{\tiny ;z),Z}%
_{1}^{(1)}{\tiny (i,z}_{0}{\tiny ;z))}}$ &  & $_{{\tiny (Z}_{2}^{(1)}%
{\tiny (1,z}_{0}{\tiny ;z),Z}_{2}^{(1)}{\tiny (i,z}_{0}{\tiny ;z))}}$ &  &  &
\\
& $_{\overset{\frac{d_{(g,i/g)}}{dz}}{\searrow}}$ &  & $_{_{\overset
{\frac{d_{(F_{1},G_{1})}}{dz}}{\searrow}}}$ &  & $_{_{\overset{\frac
{d_{(F_{2},G_{2})}}{dz}}{\searrow}}}$ & ${\tiny \vdots}$ & \\
$_{{\tiny (Z}^{(0)}{\tiny (1,z}_{0}{\tiny ;z),Z}^{(0)}{\tiny (i,z}%
_{0}{\tiny ;z))}}$ &  & $\underset{\underset{(F_{1},G_{1})}{\shortparallel}%
}{_{{\tiny (Z}_{1}^{(0)}{\tiny (1,z}_{0}{\tiny ;z),Z}_{1}^{(0)}{\tiny (i,z}%
_{0}{\tiny ;z))}}}$ &  & $\underset{\underset{(F_{2},G_{2})}{\shortparallel}%
}{_{{\tiny (Z}_{2}^{(0)}{\tiny (1,z}_{0}{\tiny ;z),Z}_{2}^{(0)}{\tiny (i,z}%
_{0}{\tiny ;z))}}}$ &  & $\underset{\underset{(F_{3},G_{3})}{\shortparallel}%
}{_{{\tiny (Z}_{3}^{(0)}{\tiny (1,z}_{0}{\tiny ;z),Z}_{3}^{(0)}{\tiny (i,z}%
_{0}{\tiny ;z))}}}$ & $_{{\tiny \ldots}}$%
\end{tabular}

In order to obtain the successor $(F_{1},G_{1})$ one can apply the
differential operator $\frac{d_{(g,i/g)}}{dz}$ to the pair of formal powers
$(Z^{(1)}(1,z_{0};z),Z^{(1)}(i,z_{0};z))$ obtaining $(Z_{1}^{(0)}%
(1,z_{0};z),Z_{1}^{(0)}(i,z_{0};z))$ which can be chosen as $(F_{1},G_{1})$.
Then this newly obtained generating pair serves for obtaining $(F_{2},G_{2})$
(differentiating $(Z_{1}^{(1)}(1,z_{0};z),Z_{1}^{(1)}(i,z_{0};z))$ in the
sense of Bers with respect to $(F_{1},G_{1})$) and positive formal powers of
subindex $2$, and in this way the whole generating sequence corresponding to
(\ref{VekuamainForg}) can be constructed.

As an illustration of the algorithm, we consider the example from the
preceding section. The generating pair $(y^{2},i/y^{2})$ was used to obtain
formal powers of order $n=0,1,2$ for $Z_{f}^{(n)}(\alpha,z_{o};z)$ of the
Vekua equation (\ref{MainVekuaf}). Then, using the transplant operator, the
corresponding formal powers $Z_{g}^{(n)}(\alpha,z_{o};z)$ of the Vekua
equation (\ref{MainVekuag}) were obtained. Looking now for a generating
sequence corresponding to the Vekua equation (\ref{MainVekuag}) we already
have $(F,G)=(g,i/g)$, where we recall that $g=\displaystyle\frac{1+xy^{3}}{y}%
$. To obtain other elements of the generating sequence for the Vekua equation
(\ref{MainVekuag}), we follow the algorithm presented above. We have%

\begin{align*}
(F_{1},G_{1}) &  =\frac{d_{(g,i/g)}}{dz}\Big(Z_{g}^{(1)}(1,z_{o}%
;z),Z_{g}^{(1)}(i,z_{o};z)\Big)\\
&  =\frac{d}{dz}\Big(Z_{g}^{(1)}(1,z_{o};z),Z_{g}^{(1)}(i,z_{o}%
;z)\Big)-A_{(g,i/g)}\Big(Z_{g}^{(1)}(1,z_{o};z),Z_{g}^{(1)}(i,z_{o}%
;z)\Big)-B_{(g,i/g)}\overline{\Big(Z_{g}^{(1)}(1,z_{o};z),Z_{g}^{(1)}%
(i,z_{o};z)\Big)}%
\end{align*}
where we used equation (\ref{FGder}) for the $(g,i/g)$-derivative in the sense
of Bers. As $A_{(g,i/g)}=0$ and
\[
B_{(g,i/g)}=\displaystyle\frac{1}{2}\displaystyle\frac{y^{4}-2ixy^{3}%
+i}{y(1+xy^{3})}%
\]
we obtain
\[
F_{1}=\displaystyle\frac{y}{y_{0}^{2}(1+xy^{3})}\Big[y(1+x_{0}y^{3}%
)-3i(x-x_{0})\Big]\quad\quad\quad\text{ and }\quad\quad\quad G_{1}%
=\displaystyle\frac{y}{3y_{0}(1+xy^{3})}\Big[y(y^{3}-y_{0}^{3})+3i(1+xy_{0}%
^{3})\Big].
\]
One can verify that $(F_{1},G_{1})$ satisfies the required property for a
generating pair on the considered domain $\Omega=\{(x,y)\ |\ x,y>0\}$:
\[
\text{Im}(\overline{F}_{1}G_{1})=\displaystyle\frac{1}{2}\left(
\displaystyle\frac{y}{y_{0}}\right)  ^{2}\displaystyle\frac{g(z_{0})}%
{g(z)}>0\text{\qquad in }\Omega.
\]

Looking now for $(F_{2},G_{2})$ we have first to calculate $Z_{g,1}%
^{(1)}(\alpha,z_{o};z)$:
\begin{align*}
\Big(Z_{g,1}^{(1)}(1,z_{o};z),Z_{g,1}^{(1)}(i,z_{o};z)\Big) &  =\frac
{d_{(g,i/g)}}{dz}\Big(Z_{g}^{(2)}(1,z_{o};z),Z_{g}^{(2)}(i,z_{o};z)\Big)\\
&  =\frac{d}{dz}\Big(Z_{g}^{(2)}(1,z_{o};z),Z_{g}^{(2)}(i,z_{o}%
;z)\Big)-B_{(g,i/g)}\overline{\Big(Z_{g}^{(2)}(1,z_{o};z),Z_{g}^{(2)}%
(i,z_{o};z)\Big)}.
\end{align*}
We obtain
\begin{align*}
Z_{g,1}^{(1)}(1,z_{o};z) &  =\displaystyle\frac{y}{15y_{0}^{2}(1+xy^{3}%
)}\Big[y\Big(15y^{3}(x^{2}-x_{0}^{2})+30(x-x_{0})+3y^{5}-5y_{0}^{2}%
y^{3}+2y_{0}^{5}\Big)\\
&  +i\Big(15(y^{2}-y_{0}^{2})+6x(y^{5}-y_{0}^{5})-45(x-x_{0})^{2}\Big)\Big]
\end{align*}
and
\[
Z_{i,g}^{(1)}(1,z_{o};z)=\displaystyle\frac{2}{3y_{0}y(1+xy^{3})}%
\Big[-(y^{3}-y_{0}^{3})(1+x_{0}y^{3})+3i(x-x_{0})y^{2}(1+xy_{0}^{3})\Big].
\]
The generating pair $(F_{2},G_{2})$ is then given by
\begin{align*}
(F_{2},G_{2}) &  =\frac{d_{(F_{1},G_{1})}}{dz}\Big(Z_{g,1}^{(1)}%
(1,z_{o};z),Z_{g,1}^{(1)}(i,z_{o};z)\Big)\\
&  =\frac{d}{dz}\Big(Z_{g,1}^{(1)}(1,z_{o};z),Z_{g,1}^{(1)}(i,z_{o}%
;z)\Big)-A_{(F_{1},G_{1})}\Big(Z_{g,1}^{(1)}(1,z_{o};z),Z_{g,1}^{(1)}%
(i,z_{o};z)\Big)-B_{(F_{1},G_{1})}\overline{\Big(Z_{g,1}^{(1)}(1,z_{o}%
;z),Z_{g,1}^{(1)}(i,z_{o};z)\Big)}%
\end{align*}
where
\[
A_{(F_{1},G_{1})}=-\frac{y^{4}+3i}{2y(1+xy^{3})}\quad\quad\quad\text{ and
}\quad\quad\quad B_{(F_{1},G_{1})}=-\frac{i}{y}.
\]
Combining these results we find
\[
F_{2}=2\left(  \frac{y}{y_{0}}\right)  ^{2}\quad\quad\quad\text{ and }%
\quad\quad\quad G_{2}=2i\left(  \frac{y_{0}}{y}\right)  ^{2},
\]
which obviously satisfies the required property that $\text{Im}(\overline
{F}_{2}G_{2})>0$. Notice that in this special case the obtained pair
$(F_{2},G_{2})$ is equivalent (in the sense introduced in \cite{Berskniga},
see also \cite{APFT}) to the  generating pair $(y^{2},i/y^{2})$ which was used
for starting the proposed procedure. As it depends on the variable $y$ only
the succeeding generating pairs can be chosen again equal to it. Thus, in this
example we constructed a complete generating sequence into which the
generating pair $(F,G)=(g,i/g)$ is embedded. Namely, $(F_{0},G_{0}%
)=(\frac{1+xy^{3}}{y},\frac{iy}{1+xy^{3}})$, $(F_{1},G_{1})=\left(  \frac
{y}{y_{0}^{2}(1+xy^{3})}\Big[y(1+x_{0}y^{3})-3i(x-x_{0})\Big],\frac{y}%
{3y_{0}(1+xy^{3})}\Big[y(y^{3}-y_{0}^{3})+3i(1+xy_{0}^{3})\Big]\right)  $,
$(F_{n},G_{n})=(y^{2},i/y^{2})$ for $n=2,3,\ldots$.

\section{Construction of Cauchy kernels\label{SectCauchykernels}}

In Section \ref{SectTransplant} we showed that the transplant operator
$T_{f,g}$ transforms a Cauchy kernel $Z_{f}^{(-1)}(\alpha,z_{0};z)$
corresponding to equation (\ref{VekuamainA}) into a Cauchy kernel
$Z_{g}^{(-1)}(\alpha,z_{0};z)$ corresponding to equation (\ref{VekuamainForg})
when $f$ and $g$ are solutions of a same Schr\"{o}dinger equation
(\ref{Schrod}). Here we give an example of the application of this procedure.

Let us consider the following Vekua equation
\begin{equation}
w_{\overline{z}}=\frac{1}{2}\left(  \frac{1}{x}+\frac{i}{y}\right)
\overline{w} \label{VekuaCauchy}%
\end{equation}
in the domain $\Omega=\{(x,y)\ |\ x>0\mbox{ and }y>0\}$. Note that the
coefficient in the equation admits a representation in the form of a
logarithmic derivative of a real-valued function:
\[
\frac{1}{2}\left(  \frac{1}{x}+\frac{i}{y}\right)  =\frac{g_{\overline{z}}}{g}%
\]
where moreover, $g=xy$ is a harmonic function. The simplest nontrivial
harmonic function is, of course, $f\equiv1$. The corresponding Vekua equation
(\ref{VekuamainA}) is just the Cauchy-Riemann system for which the Cauchy
kernel is well known. In order to obtain $Z_{g}^{(-1)}(\alpha,z_{0};z)$ for
any $\alpha\in\mathbb{C}$ we need to calculate $Z_{g}^{(-1)}(1,z_{0};z)$ and
$Z_{g}^{(-1)}(i,z_{0};z)$ applying the transplant operator $T_{1,xy}$ to
$1/(z-z_{0})$ and $i/(z-z_{0})$, respectively. We will show the result for
$Z_{g}^{(-1)}(1,z_{0};z)$ (an expression for $Z_{g}^{(-1)}(i,z_{0};z)$ can be
obtained analogously) calculated with the aid of Maple.

Thus we consider $W=1/(z-z_{0})$ for $(x_{0},y_{0})\in\Omega$ and find
\[
Z_{g}^{(-1)}(1,z_{0};z)=T_{1,xy}[W]=\frac{x-x_{0}}{|z-z_{0}|^{2}}+\frac{i}%
{xy}\overline{A}\Big[i(xy)^{2}\partial_{\overline{z}}\Big(\frac{x-x_{0}%
}{xy|z-z_{0}|^{2}}\Big)\Big],
\]
where
\begin{align*}
i(xy)^{2}\partial_{\overline{z}}\Big(\frac{x-x_{0}}{xy|z-z_{0}|^{2}}\Big)  &
=\frac{1}{2|z-z_{0}|^{4}}\Big[\Big(y_{0}^{2}x^{2}-x_{0}^{3}x+3x^{2}%
y^{2}-3x_{0}x^{3}+3x_{0}^{2}x^{2}+4x_{0}y_{0}xy\\
&  -4y_{0}x^{2}y-3x_{0}xy^{2}-x_{0}y_{0}^{2}x+x^{4}\Big)-i\Big(-x_{0}y_{0}%
^{2}y+4x_{0}^{2}xy\\
&  +2x_{0}y_{0}y^{2}-5x_{0}x^{2}y+2x^{3}y-x_{0}^{3}y-x_{0}y^{3}\Big)\Big].
\end{align*}
Applying the $\overline{A}$ operator (integrating from $(0,0)$ to $(x,y)$) we
obtain
%

\begin{align*}
\overline{A}\Big[i(xy)^{2}\partial_{\overline{z}}\Big(\frac{x-x_{0}%
}{xy|z-z_{0}|^{2}}\Big)\Big]  &  =\frac{1}{2|z-z_{0}|^{2}}\Big\{\big(2x_{0}%
^{2}y_{0}-4y_{0}^{2}y-4x_{0}y_{0}x+2y_{0}y^{2}+2y_{0}^{3}+2y_{0}%
x^{2}\big)\arctan(\frac{y_{0}}{x_{0}})\\
&  +\big(-2x_{0}^{2}x+x_{0}x^{2}+x_{0}^{3}+x_{0}y^{2}+x_{0}y_{0}^{2}%
-2x_{0}y_{0}y\big)\ln\left\vert \frac{z-z_{0}}{z}\right\vert ^{2}\\
&  +\big(-2x_{0}^{2}y_{0}+4y_{0}^{2}y+4x_{0}y_{0}x-2y_{0}y^{2}-2y_{0}%
^{3}-2y_{0}x^{2}\big)\arctan(\frac{y-y_{0}}{x-x_{0}})\\
&  +\big(-2y_{0}xy+2x_{0}^{2}x+2y_{0}^{2}x+2x^{3}-4x_{0}x^{2}\big)\Big\}+c.
\end{align*}
Choosing $c=0$ we find
\begin{align}
Z_{g}^{(-1)}(1,z_{0};z)  &  =T_{1,xy}[W]\nonumber\\
&  =\frac{x-x_{0}}{|z-z_{0}|^{2}}+\frac{i}{2xy|z-z_{0}|^{2}}\Big\{\big(2x_{0}%
^{2}y_{0}-4y_{0}^{2}y-4x_{0}y_{0}x+2y_{0}y^{2}+2y_{0}^{3}+2y_{0}%
x^{2}\big)\arctan(\frac{y_{0}}{x_{0}})\nonumber\\
&  +\big(-2x_{0}^{2}x+x_{0}x^{2}+x_{0}^{3}+x_{0}y^{2}+x_{0}y_{0}^{2}%
-2x_{0}y_{0}y\big)\ln\left\vert \frac{z-z_{0}}{z}\right\vert ^{2}\nonumber\\
&  +\big(-2x_{0}^{2}y_{0}+4y_{0}^{2}y+4x_{0}y_{0}x-2y_{0}y^{2}-2y_{0}%
^{3}-2y_{0}x^{2}\big)\arctan(\frac{y-y_{0}}{x-x_{0}})\nonumber\\
&  +\big(-2y_{0}xy+2x_{0}^{2}x+2y_{0}^{2}x+2x^{3}-4x_{0}x^{2}\big)\Big\}.
\label{Zminus1}%
\end{align}

\begin{figure}[H]
\centering
\includegraphics[width=10cm]{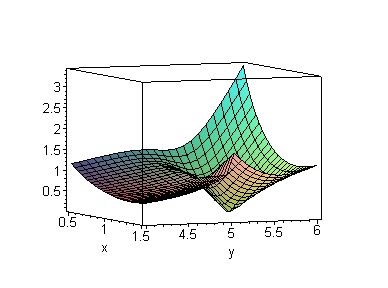}
\caption{Graph of the function $H(x,y)=\big|(z-z_0)\,Z_g^{(-1)}(1,z_0;z)\big|$ for the function $Z_g^{(-1)}(1,z_0;z)$ given by (\ref{Zminus1}) and $(x_0,y_0)=(1,5)$.}
\label{fig}
\end{figure}

Property (\ref{condCauchykernel}) for $Z_{g}^{(-1)}(1,z_{0};z)$ is illustrated
in Figure~\ref{fig} by the graph of the function $H(x,y)=\left| \frac
{Z_{g}^{(-1)}(1,z_{0};z)}{(z-z_{0})^{-1}}\right| $ for $(x_{0},y_{0}%
)=(1,5)\in\Omega$.

\section{Hyperbolic pseudoanalytic function theory}

In \cite{KravRocTre} and \cite{KKT} \textquotedblleft hyperbolic
pseudoanalytic function theory\textquotedblright\ was studied where hyperbolic
numbers $\mathbf{D}$ (also called duplex numbers) \cite{Sobczyk} defined by
\[
\mathbf{D}:=\big\{x+tj\ :\ j^{2}=1,\ x,t\in\mathbb{R}\big\}\cong%
\mathrm{Cl}_{\mathbb{R}}(0,1)
\]
are considered instead of (elliptic) complex numbers. Here we show that the
concept of the transplant operator can be introduced in this context as well
and it allows one to solve hyperbolic main Vekua equations related to the
Klein-Gordon equation.

As in the case of complex numbers, we denote the real and imaginary parts of
$z=x+tj\in\mathbf{D}$ by $x=\mathrm{Re}~z$ and $t=\mathrm{Im}~z$. Now, by
defining the conjugate as $\bar{z}=Cz:=x-tj$ and the hyperbolic modulus as
$|z|^{2}:=z\bar{z}=x^{2}-t^{2}$, we can verify that the inverse of $z$
whenever exists is given by
\[
z^{-1}=\displaystyle\frac{\overline{z}}{|z|^{2}}.
\]
The set $\mathcal{NC}$ of zero divisors for hyperbolic numbers $\mathbf{D}$,
called the \emph{null-cone}, is given by $\mathcal{NC}%
=\big\{x+tj\ :\ |x|=|t|\big\}$.

\begin{definition}
Let $U$ be an open set in $\mathbf{D}$ and $z_{0}\in U$. Then $w:U\subseteq
\mathbf{D}\longrightarrow\mathbf{D}$ is said to be $\mathbf{D}$-differentiable
at $z_{0}$ with derivative equal to $w^{\prime}(z_{0})\in\mathbf{D}$ if
\[
\lim_{\overset{\scriptstyle z\rightarrow z_{0}}{\scriptscriptstyle(z-z_{0}%
\mbox{
}inv.)}}\frac{w(z)-w(z_{0})}{z-z_{0}}=w^{\prime}(z_{0}).
\]

\end{definition}

Here $z$ tends to $z_{0}$ following the invertible trajectories. We also say
that the function $w$ is $\mathbf{D}$-holomorphic on an open set $U$ if and
only if $w$ is $\mathbf{D}$-differentiable at each point of $U$.

\begin{theorem}
Let $U$ be an open set and $w:U\subseteq\mathbf{D}\longrightarrow\mathbf{D}$
such that $w\in{C}^{1}(U)$. Let also $w(x+tj)=w_{1}(x,t)+w_{2}(x,t)j$. Then
$w$ is $\mathbf{D}$-holomorphic on $U$ if and only if%

\begin{equation}
\frac{\partial{w_{1}}}{\partial{x}}=\frac{\partial{w_{2}}}{\partial{t}%
}\ \ \ \ \ \mbox{ and }\ \ \ \ \ \frac{\partial{w_{2}}}{\partial{x}}%
=\frac{\partial{w_{1}}}{\partial{t}}. \label{hyperbolicC-R}%
\end{equation}
Moreover $w^{\prime}=\displaystyle\frac{\partial{w_{1}}}{\partial{x}%
}+\displaystyle\frac{\partial{w_{2}}}{\partial{x}}j$ and $w^{\prime}(z)$ is
invertible if and only if $\det\mathcal{J}_{w}(z)\neq0$, where $\mathcal{J}%
_{w}(z)$ is the Jacobian matrix of $w$ at $z$. \label{theobasic}
\end{theorem}

System of equations (\ref{hyperbolicC-R}) is called \textquotedblleft
hyperbolic Cauchy-Riemann\textquotedblright\ equations. It was considered in
\cite{18}, \cite{31}, \cite{32}.

For $z=x+tj$ where $x,t$ are real variables, we define the operators
$\partial_{z}$ and $\partial_{\overline{z}}$ in the hyperbolic function theory
as
\[
\partial_{z}=\frac{1}{2}\left(  {\partial_{x}+j\partial_{t}}\right)
\quad\quad\mbox{ and }\quad\quad\partial_{\bar{z}}=\frac{1}{2}\left(
{\partial_{x}-j\partial_{t}}\right)  ,
\]
such that $w_{\overline{z}}(z)=0$ if and only if hyperbolic Cauchy-Riemann
equations (\ref{hyperbolicC-R}) are satisfied.

Let $\Omega$ be a domain in $\mathbf{R}^{2}$ without zero divisors. We
consider now the hyperbolic Vekua equation
\begin{equation}
W_{\overline{z}}=\frac{f_{\overline{z}}}{f}\overline{W}\text{\qquad in }%
\Omega, \label{hypVekua}%
\end{equation}
where $f$ is a positive function of $x$ and $t$, twice continuously
differentiable, which will be supposed to be a particular solution of the
following $(1+1)$-dimensional Klein-Gordon equation
\begin{equation}
(\square-q)\varphi=0\text{\qquad in }\Omega. \label{K-G}%
\end{equation}
Here $\square:=\displaystyle\frac{\partial^{2}}{\partial x^{2}}-\frac
{\partial^{2}}{\partial t^{2}}$, the potential $q$ is a real valued function
and $\varphi$ is a twice continuously differentiable real valued function of
$x$ and $t$.

\begin{theorem}
\cite{KravRocTre} \label{thmfactoKlein} Let $f$ be a positive particular
solution of (\ref{K-G}) in $\Omega$. Then for any real valued function
$\varphi\in C^{2}(\Omega)$ the following equalities hold
\[
\frac{1}{4}(\square-q )\varphi= \big(\partial_{\bar{z}}+\displaystyle\frac{
f_{z}}{f}C\big)\big(\partial_{z}-\displaystyle\frac{f_{z}}{f}C\big)\varphi=
\big(\partial_{z}+\displaystyle\frac{f_{\bar{z}}}{f}C\big)\big(
\partial_{\bar{z}}-\displaystyle\frac{f_{\bar{z}}}{f}C\big)\varphi.
\label{factoKlein}
\]

\end{theorem}

In a similar way as in the elliptic case, an immediate consequence of this
theorem is the fact that if $W$ is a solution of (\ref{hypVekua}) then its
real part $W_{1}$ is a solution of (\ref{K-G}), meanwhile its imaginary part
$W_{2}$ is a solution of the following Klein-Gordon equation
\begin{equation}
(\square-q_{1})W_{2}=0\qquad\text{in }\Omega, \label{K-G2}%
\end{equation}
where $q_{1}=-q+8\displaystyle\frac{|f_{z}|^{2}}{f^{2}}$ (see\cite{KravRocTre}%
, \cite{APFT}).

Note that in the hyperbolic case the family of real valued functions $\phi$
such that $\partial_{\bar{z}}\phi=\Phi$, and $\Phi=\Phi_{1}+j\Phi_{2}$, can be
constructed as
\[
\overline{A}_{h}[\Phi](x,t)=2\left(  \int_{x_{0}}^{x}\Phi_{1}(\eta
,t)\mathrm{d}\eta-\int_{t_{0}}^{t}\Phi_{2}(x_{0},\xi)\mathrm{d}\xi\right)
+c,
\]
when $\Phi$ satisfies the compatibility condition
\[
\partial_{t}\Phi_{1}+\partial_{x}\Phi_{2}=0.
\]

\begin{theorem}
\cite{KravRocTre} Given a solution $W_{1}$ of the Klein-Gordon equation
(\ref{K-G}), the corresponding $W_{2}$ such that $W=W_{1}+jW_{2}$ is a
solution of the hyperbolic Vekua equation (\ref{hypVekua}) can be constructed
according to the formula
\[
W_{2}=-f^{-1}\overline{A}_{h}\big[jf^{2}\partial_{\bar{z}}(f^{-1}W_{1})\big].
\]

Vice versa, given a solution $W_{2}$ of the Klein-Gordon equation
(\ref{K-G2}), the corresponding $W_{1}$ such that $W=W_{1}+jW_{2}$ is a
solution of the hyperbolic Vekua equation (\ref{hypVekua}) can be constructed
according to the formula
\[
W_{1}=-f\overline{A}_{h}\big[jf^{-2}\partial_{\bar{z}}(fW_{2})\big].
\]

\end{theorem}

Now, let $g$ be another positive solution of (\ref{K-G}) associated with the
following hyperbolic Vekua equation
\begin{equation}
w_{\overline{z}}=\frac{g_{\overline{z}}}{g}\overline{w}\text{\qquad in }%
\Omega. \label{hypVekua2}%
\end{equation}
For $W$ and $w$ solutions of hyperbolic Vekua equations (\ref{hypVekua}) and
(\ref{hypVekua2}), respectively, we have that both $\mbox{Re}\,W$ and
$\mbox{Re}\,w$ satisfy (\ref{K-G}), meanwhile $\mbox{Im}\,W$ and
$\mbox{Im}\,w$ satisfy two different Klein-Gordon equations
\begin{align*}
(\square-q_{1})\mbox{Im}\,W  &  =0\qquad\text{in }\Omega\\
(\square-q_{2})\mbox{Im}\,w  &  =0\qquad\text{in }\Omega,
\end{align*}
respectively, where $q_{1}=-q+8\displaystyle\frac{|f_{z}|^{2}}{f^{2}}$ and
$q_{2}=-q+8\displaystyle\frac{|g_{z}|^{2}}{g^{2}}$.

In a similar way as in the elliptic case, we introduce a hyperbolic transplant
operator $\widetilde{T}_{f,g}$ which transforms solutions of (\ref{hypVekua})
into solutions of (\ref{hypVekua2}) in the following way:
\[
\widetilde{T}_{f,g}[W]=P^{+}W-jg^{-1}\overline{A}_{h}[jg^{2}\partial
_{\overline{z}}(g^{-1}P^{+}W)]
\]
where $P^{+}=\frac{1}{2}(I+C)$.

Again by assigning a fixed value in a certain point of $\Omega$ to the result
of application of $\overline{A}_{h}$, we obtain an invertible one-to-one map
establishing a relation between solutions of (\ref{hypVekua}) and
(\ref{hypVekua2}). The inverse of $\widetilde{T}_{f,g}$ is given by the
expression
\[
\widetilde{T}_{f,g}^{-1}[w]=\widetilde{T}_{g,f}[w]=P^{+}w-jf^{-1}\overline
{A}_{h}[jf^{2}\partial_{\overline{z}} (f^{-1}P^{+}w)].
\]

The nonnegative formal powers $Z_{m}^{(n)}(a,z_{0};z)$, where $a,z_{0}$ and
$z$ hyperbolic numbers, are defined in hyperbolic pseudoanalytic function
theory as in definition \ref{DefFormalPower} for usual (elliptic)
pseudoanalytic theory. These formal powers have same properties replacing $i$
by $j$ everywhere \cite{KravRocTre}. Therefore, using hyperbolic transplant
operator as in the elliptic case, nonnegative formal powers and generating
sequence of the hyperbolic Vekua equation (\ref{hypVekua2}) can be obtained
from nonnegative formal powers and generating sequence of (\ref{hypVekua}) in
a domain $\Omega$.

\end{document}